\documentclass[11pt]{article}

\usepackage{amssymb}

\newtheorem{theorem}{Theorem}

\newtheorem{lemma}[theorem]{Lemma}

\newtheorem{prop}[theorem]{Proposition}
\newtheorem{remark}{Remark}

\def\dde{\stackrel{{\rm def}}{=}}
\def\elaw{\stackrel{d}{=}}

\def\AA{{\bf A}}
\def\AAA{{\cal A}}
\def\ee{\varepsilon}
\def\esp{{\mathbb{E}}}
\def\FF{{\cal F}}
\def\HH{{\cal H}}
\def\k{{\cal K}}
\def\lacc{\left\{}
\def\lcr{\left[}
\def\LL{{\cal L}}
\def\lpa{\left(}
\def\lva{\left|}
\def\mude{{\tilde \mu}}
\def\NN{{\mathbb{N}}}
\def\pb{{\mathbb{P}}}
\def\R{{\mathbb{R}}}
\def\racc{\right\}}
\def\rcr{\right]}
\def\rpa{\right)}
\def\rva{\right|}
\def\Un{{{1\mskip-4mu {\rm l}}}}
\def\TT{{\bf T}}
\def\VV{{\cal V}}
\def\WW{{\cal W}}

\newcommand{\fin}{\vspace{-0.3cm}
                  \begin{flushright}
                  \mbox{$\Box$}
                  \end{flushright}
                  \noindent}

\begin{document}

\begin{center}

\huge
{\bf Small ball estimates in $p$-variation for stable processes}

\vspace{8mm}

\Large
{\bf Thomas Simon}

\vspace{5mm}
\normalsize
{\em 
Equipe d'Analyse et Probabilit\'es, Universit\'e d'Evry-Val d'Essonne\\
Boulevard Fran\c{c}ois Mitterrand, F-91025 EVRY cedex\\
e-mail:} {\tt simon@maths.univ-evry.fr}   
\end{center}

\vspace{5mm}

\begin{abstract}

\vspace{2mm}

\noindent
{Let $\{Z_t, \; t\geq 0\}$ be a strictly stable process on $\R$ with index $\alpha\in (0,2]$. We prove that for every $p > \alpha$, there exists $\gamma = \gamma (\alpha, p)$ and $\k = \k (\alpha, p)\in (0, +\infty)$ such that 
$$\lim_{\ee\downarrow 0}\ee^{\gamma}\log\pb\lcr ||Z||_{p}\leq \ee \rcr
 \; =\; - \k,$$
where $||Z||_{p}$ stands for the strong $p$-variation of $Z$ on $[0,1]$. The critical exponent $\gamma (\alpha, p)$ takes a different shape according as $|Z|$ is a subordinator and $p >1$, or not. The small ball constant $\k (\alpha, p)$ is explicitly computed when $p \leq 1$, and a lower bound on $\k (\alpha, p)$ is easily obtained in the general case. In the symmetric case and when $p > 2$, we can also give an upper bound on $\k (\alpha, p)$ in terms of the Brownian small ball constant under the $(1/p)$-H\"older semi-norm. Along the way, we remark that the positive random variable $||Z||^p_{p}$ is not necessarily stable when $p > 1$, which gives a negative answer to an old question of P.~E.~Greenwood \cite{Gr}.}
\end{abstract}

\vspace{4mm}

\noindent
{\bf Keywords:} H\"older semi-norm - $p$-variation - Small balls probabilities - Stable processes - Subordination.

\vspace{2mm}

\noindent
{\bf MSC 2000:} 60F99, 60G52

\section{Introduction}

In a recent paper \cite{LiSi}, a general method was introduced to prove the existence of finite small ball constants for real fractional $\alpha$-stable processes, under different norms. Whereas the existence result holds with a reasonable level of generality including for example all symmetric $\alpha$-stable processes - see Theorem 3.1 in \cite{LiSi}, the finiteness result (which amounts to a lower estimate on the small ball probabilities - see Theorem 4.1 in \cite{LiSi}) was obtained with the help of wavelet decompositions concerning only {\em continuous} processes. The suitable lower estimate for small probabilities of $\alpha$-stable L\'evy processes under the uniform norm is a classical result, which dates back to Taylor \cite{Ta} and Mogul'ski\v {\i} \cite{Mo}. By comparison, this lower estimate entails immediately the good lower bounds under all the $\LL_p$-norms ($1\leq p\leq \infty$), since the critical exponent $\gamma = \alpha$ does not depend on $p$.

A natural (and finer) semi-norm on the set of real c\`ad-l\`ag functions is the strong $p$-variation in the sense of N. Wiener, which was quite intensively studied by the stochastic community in the midst of the last century. Bretagnolle \cite{Bre} had obtained a general criterion ensuring that an $\alpha$-stable L\'evy process has a.s. finite $p$-variation if and only if $p > \alpha$ (see also \cite{Le} and \cite{BG} for previous results in the symmetric case). More recently, Chistyakov and Galkin \cite{ChG} proved an interesting embedding theorem which entails that for continuous paths, $p$-variation and $(1/p)$-H\"older semi-norm are roughly equivalent notions when $p\geq 1$. In the discontinuous framework, $p$-variation seems to be a good substitute for the irrelevant $(1/p)$-H\"older semi-norm, in studying finer sample path properties. We finally refer to the comprehensive survey of Dudley and Norvai\v sa \cite{DN} for recent developments about this notion, both from the analytical and probabilistic point of view. 

The purpose of the present paper is to prove the existence of the small deviation constant for strictly (non necessarily symmetric) $\alpha$-stable processes with respect to the $p$-variation $||.||_p$, when $p > \alpha$. In other words, we prove that 
\begin{equation}
\label{smdev}
\lim_{\ee\downarrow 0}\ee^{\gamma_{\alpha,p}}\log\pb\lcr ||Z||_{p}\leq \ee \rcr\; =\; - \k (\alpha, p)\,\in\, (0, +\infty),
\end{equation}
where the critical exponent is given by $\gamma_{\alpha,p} = p\alpha/(p-\alpha)$, except when $|Z|$ is a subordinator and $p >1$, where $\gamma_{\alpha,p} = \alpha/(1-\alpha)$. This confirms the prediction of Section 6.3 in \cite{LiSi}.

Notice that there are many ways to define the variation of a function, and for $\alpha$-stable processes this question had been thoroughly studied by Greenwood \cite{Gr} and Fristedt-Taylor \cite{FT} (see also the references therein). Greenwood had pointed out that when $\alpha < p\leq 1$, the strong $p$-variation of an $\alpha$-stable process is the $(1/p)$-th power of some positive $(\alpha/p)$-stable variable. Actually, it is easy to determine the scaling parameter of this latter variable, and then De Bruijn's exponential Tauberian theorem yields (\ref{smdev}) readily, with an explicit formula for the small ball constant $\k (\alpha, p)$ as a bonus. The same method holds for the case when $|Z|$ is a subordinator and $p >1$, because here $||Z||_p = |Z_1|$.

In all the other situations, we first prove that the random variable $||Z||^p_p$ is not stable. This gives a negative answer to a question of Greenwood - see Section 6 in \cite{Gr}, and makes the solution to our small deviation problem more involved. The proof of (\ref{smdev}) is carried out with classical scaling arguments relying on an appropriate discretization similar to those in \cite{Ta} and \cite{Mo}, and on the rough positivity results recently obtained in \cite{pvaria}. 

When $|Z|$ is not a subordinator and $p >1$, we were not able to compute the small ball constant $\k (\alpha, p)$, even when $Z$ is Brownian motion. When $\alpha <2$, a general positive lower bound is easily obtained by comparison with the sum of the $p$-th power of the jumps. Thanks to the equality of the critical exponents, the upper estimate $\k (2, p) \leq c_p$ follows readily, where $c_p$ is the small constant for Brownian motion under the $(1/p)$-H\"older semi-norm. However, no positive lower bound is available for $\k (2, p)$. In the symmetric case and when $\alpha < 2 < p$, an algebraic surprise arising from Bochner's subordination makes it possible to obtain an upper bound on $\k (\alpha, p)$ in terms of $c_p$ and the (explicit) constant for $(\alpha/2)$-stable subordinators. 

We remark that the above arguments also yield two-sided estimates on the small ball constant under the supremum and oscillation semi-norms in the symmetric Non-Gaussian framework. One may ask if these arguments could be refined and yield sharper results, perhaps with the help of more advanced stochastic calculus as in \cite{Shi}. We prefer leaving this question open for now. 

\section{Preliminaries and statement of the results}

\subsection{Strictly stable processes and their strong $p$-variation}

We consider $\lacc Z_t, \; t\geq 0\racc$ a strictly stable process on $\R$ with index $\alpha\in(0,2]$, viz. a real L\'evy process satisfying the following self-similarity property: for every $k > 0$
$$\lacc Z_{kt}, \; t\geq 0\racc\;\elaw\;\lacc k^{1/\alpha}Z_t, \; t\geq 0\racc.$$
When $\alpha = 2$, $Z$ is just a rescaled Brownian motion: there exists $a\in\R$ such that $Z = a W,$ where $\lacc W_t, \; t\geq 0\racc$ is a standard linear Brownian motion. When $\alpha < 2$, $Z$ has the following L\'evy-It\^o decomposition: there exists $b\in\R$, $c_-, c_+\in\R^+$ such that
\begin{equation}
\label{li}
Z_t\; = \; b t \; + \; \int_0^t\int_{|z| \leq 1}\!\!
z\,\mude(ds,dz) \; + \; \int_0^t\int_{|z| > 1}\!\! z\,\mu(ds,dz)
\end{equation}
for every $t \geq 0$, where $\nu$ is a measure on $\R^*$ given by
$$\nu(dz)\; =\; \lacc\begin{array}{ll} c_+z^{-\alpha -1} dz & \mbox{if $z > 0$}\\
                                       c_-|z|^{-\alpha -1} dz & \mbox{if $z < 0$}
                     \end{array}
                \right.,
$$ 
$\mu$ is the Poisson measure over $\R^+\!\times\R^d$ with intensity $ds \otimes\nu(dz)$, and $\mude = \mu - ds \otimes\nu$ is the compensated measure. When $\alpha = 1$, the drift $b$ can take any real value. When $\alpha \neq 1$, the strict stability imposes a fixed value on $b$ given by the relation:
$$c_+\; =\; c_- + b(1-\alpha).$$
In particular, we see that when $\alpha < 1$, $Z$ is a pure jump process, i.e. it can be rewritten as the sum of its jumps:
\begin{equation}
\label{jump}
Z_t \; = \; \sum_{s\leq t} \Delta Z_s
\end{equation}
for every $t \geq 0$. In the remainder of this paper we will exclude implicitly the trivial case when $Z$ is a pure drift, i.e. we will suppose that $a\neq 0$ (resp. $\nu\not \equiv 0$) when $\alpha = 2$ (resp. $\alpha\neq 2$). For the sake of concision, ``stable'' will always mean ``strictly stable'' further on. We will always consider a c\`ad-l\`ag (resp. a continuous) version of $Z$ when $\alpha < 2$ (resp. $\alpha = 2$). Sometimes, we will focus on two particular cases:

\vspace{2mm}

\noindent
(i) The case when $Z$ is symmetric. Its Fourier transform is then given by
\begin{equation} \label{kappa}
\esp\lcr e^{\rm{i} \lambda Z_t}\rcr\; =\; e^{-\kappa t|\lambda|^{\alpha}}
\end{equation}
for every $t\geq 0$ and $\lambda\in\R$, with $\kappa > 0$ a normalization constant. When $\alpha = 2$, then $\kappa = a^2/2$ with the above notation. When $\alpha < 2$, $Z$ can be viewed as a pure jump process like in (\ref{jump}) - even here though the series on the right hand side does not converge absolutely when $\alpha \geq 1$. Besides, we have
\begin{equation} \label{cplus}
c_- = c_+ = \lpa \frac{\alpha \Gamma(\alpha) \sin (\pi\alpha/2)}{\pi}\rpa\kappa
\end{equation} 
in (\ref{li}), where $\Gamma$ is the usual Gamma function (see e.g. Lemma 14.11 in \cite{sato}).

\vspace{2mm}

\noindent
(ii) The case when $Z$ is a subordinator. Necessarily $\alpha < 1$ and the Laplace transform of $Z$ is given by 
$$\esp\lcr e^{-\lambda Z_t}\rcr\; =\; e^{-\kappa t\lambda^{\alpha}}$$
for every $t\geq 0$ and $\lambda\in\R^+$, with $\kappa > 0$ a normalization constant. Again, $Z$ can be viewed as a pure jump process, and in (\ref{li}) one has $a = c_- = 0$, whereas
\begin{equation} \label{cplus2}
c_+ = \lpa \frac{\alpha}{\Gamma(1-\alpha)}\rpa\kappa
\end{equation}
(see Example 24.12 in \cite{sato}). In particular, we see that 
\begin{equation}
\label{laplace1}
\lambda^{-\alpha}\log \esp\lcr e^{-\lambda Z_1}\rcr\; =\; - \lpa\frac{c_+\Gamma(1-\alpha)}{\alpha}\rpa
\end{equation}
for every $\lambda > 0$.

\vspace{2mm}

When $\alpha < 2$, we will make a repeated use of the following pure jump process, which takes finite values if and only if $p > \alpha$:
$$S^p_t\; =\; \sum_{s\leq t} {\lva \Delta Z_s\rva}^p$$ 
for every $t\geq 0$. When $p > \alpha$, it is easy to see (or follows from Proposition 3 in \cite{FT}) that $S^p$ is a stable subordinator with index $\alpha/p$, whose jumping measure is given by
$$\nu(dz)\; =\; \frac{(c_- + c_+)dz}{p z^{(1 + \alpha/p)}}$$
on $\R^+$. In particular, we see that 
\begin{equation}
\label{laplace2}
\lambda^{-\alpha/p}\log \esp\lcr e^{-\lambda S^p_1}\rcr\; =\; - \lpa\frac{(c_- + c_+)}{\alpha}\rpa \Gamma(1- \alpha/p)
\end{equation}
for every $\lambda > 0$.

\vspace{2mm}

For every $p > 0$, the (strong) $p$-variation of a regulated function $f : \R^+ \rightarrow \R$ over a closed finite interval $I$ of $\R^+$ is defined by
$${\lva\lva f\rva\rva}_{I, p}\; =\; {\lpa\sup_{\tiny{t_0 < \ldots < t_k\in 
I}}\sum_{j=1}^k{\lva f(t_j) - f(t_{j-1})\rva}^p\rpa}^{1/p}.$$
Thanks to Minkowski's inequality, we see that ${\lva\lva .\rva\rva}_{I, p}$ is a semi-norm when $p\geq 1$, but this is no more true when $p <1$. Actually, when $p <1$, ${\lva\lva f\rva\rva}_{I, p} = +\infty$ unless $f$ is a pure jump function, where we have
\begin{equation}
\label{quasi}
{\lva\lva f\rva\rva}_{I, p}\; =\; {\lpa \sum_{t\in I} {\lva \Delta f_t\rva}^p\rpa}^{1/p}.
\end{equation}
For the sake of concision, we will note ${\lva\lva .\rva\rva}_p = {\lva\lva .\rva\rva}_{[0,1], p}$ subsequently. From (\ref{quasi}) we see immediately that if $p < 1$, then
$${\lva\lva Z \rva\rva}_p < +\infty \;\;\mbox{a.s.}\;\Longleftrightarrow\; \alpha < p.$$
It follows from the classical results of L\'evy \cite{Le} for Brownian motion and Bretagnolle \cite{Bre} for general L\'evy processes, that the above equivalence remains true without any restriction on $p > 0$. When $\alpha = 2$, the if part is obtained immediately from the obvious inequality
\begin{equation}
\label{hoelder}
{\lva\lva f \rva\rva}_p \;\leq\; \sup_{0\leq s< t \leq 1}\frac{|f(t)-f(s)|}{{|t-s|}^{1/p}},
\end{equation}
which is valid for every function $f : [0,1]\rightarrow \R$ and every $p \geq 1$, and the fact that Brownian motion has a.s. $(1/p)$-H\"older paths as soon as $p > 2$. Similarly, one can show that Brownian motion has a.s. finite $p$-variation only if $p \geq 2$, but the exclusion of the boundary case $p =2$ requires a more subtle analysis. When $\alpha < 2$, the only if part follows easily from the straightforward inequality:
\begin{equation}
\label{infvar}
{\lva\lva f\rva\rva}_{p}\; \geq\; {\lpa \sum_{t\leq 1} {\lva \Delta f_t\rva}^p\rpa}^{1/p},
\end{equation}
which is valid for every function $f : [0,1]\rightarrow \R$ and every $p \geq 0$, and the fact that $S^p_1 = +\infty$ a.s. when $p\geq 1$. In the symmetric case, the if part is not difficult to prove via (\ref{quasi}), (\ref{hoelder}), and Bochner's subordination \cite{BG}.

\subsection{The case when $|Z|$ is a subordinator or $p\leq 1$.}

In general, little is known about the distributional properties of ${\lva\lva Z \rva\rva}_p$, even when $Z$ is Brownian motion. However, the law of ${\lva\lva Z \rva\rva}_p$ is completely explicit when $|Z|$ is a subordinator or when $\alpha < p\leq 1$. Indeed, we see from (\ref{quasi}) that if $\alpha < p < 1$, then
$${\lva\lva Z\rva\rva}^p_p\; =\; \sum_{t\geq 1} {\lva \Delta Z_t\rva}^p\; =\; S^p_1.$$
If $p\geq 1$ and $|Z|$ is a subordinator, then the monotonicity of $Z$ and the fact that $(a + b)^p \geq a^p + b^p$ for every $a, b \geq 0$ implies that
\begin{equation}\label{pge1}
{\lva\lva Z\rva\rva}^p_p\; =\; {\lva Z_1 - Z_0\rva}^p\; =\; |Z_1|^p.
\end{equation}
Finally, if $\alpha < p = 1$ and $|Z|$ is not a subordinator, then we can write $Z = Z^1 - Z^2$, where $Z^1$ and $Z^2$ are two independent stable subordinators of index $\alpha$, so that by triangle's inequality,
$${\lva\lva Z\rva\rva}_1\; \leq \; {\lva\lva Z^1\rva\rva}_1 \; + \;{\lva\lva Z^2\rva\rva}_1 \; =\; \sum_{t\geq 1} \lva \Delta Z^1_t\rva\; +\; \sum_{t\geq 1} \lva \Delta Z^2_t\rva\; =\; \sum_{t\geq 1} \lva \Delta Z_t\rva\; =\; S^1_1,$$
and it follows from (\ref{infvar}) that
$${\lva\lva Z\rva\rva}_1\; =\; \sum_{t\geq 1} \lva \Delta Z_t\rva\; =\; S^1_1.$$
In particular, we see that when $p\leq 1$, the $p$-variation of a stable process with index $\alpha < p$ is the $(1/p)$-th power of a stable positive random variable with index $\alpha/p$. This had been proved a long time ago by Greenwood - see Theorem 1 in \cite{Gr}, with different arguments. 

\vspace{2mm}

The above considerations yield easily a complete answer to our small deviation problem when $|Z|$ is a subordinator or $p\leq 1$:

\begin{prop} \label{p21}
Let $Z$ be a real stable process with index $\alpha < 1$. Then, for every $p\in (\alpha, 1]$,
\begin{equation}
\label{ap1}
\lim_{\ee\downarrow 0}\ee^{\frac{p\alpha}{p-\alpha}}\log\pb\lcr ||Z||_{p}\leq \ee \rcr
 \; =\; - \lpa \frac{p-\alpha}{\alpha}\rpa {\lpa \lpa\frac{(c_- + c_+)}{p}\rpa\Gamma(1-\alpha/p)\rpa}^{\frac{p}{p-\alpha}}
\end{equation}
with the above notations. Besides, if $|Z|$ is a subordinator, then for every $p > 1$,
\begin{equation}
\label{subo}
\lim_{\ee\downarrow 0}\ee^{\frac{\alpha}{1-\alpha}}\log\pb\lcr ||Z||_{p}\leq \ee \rcr
 \; =\; - \lpa 1/\alpha - 1\rpa {\lpa c_{\pm} \Gamma(1-\alpha)\rpa}^{\frac{1}{1-\alpha}}
\end{equation}
with the obvious notation for $c_{\pm}$.
\end{prop}

\noindent
{\em Proof.} When $p\in (\alpha, 1]$, we saw that ${\lva\lva Z\rva\rva}^p_p\; =\; S^p_1,$ so that (\ref{ap1}) follows from (\ref{laplace2}) and De Bruijn's exponential Tauberian theorem (see Theorem 3.5 in \cite{LS}, or Theorem 4.12.9 in \cite{BGT} for the most general formulation). When $|Z|$ is a subordinator, then ${\lva\lva Z\rva\rva}_p\; =\; Z_1$ for every $p > 1$, so that (\ref{subo}) follows from (\ref{laplace1}) in the same way.

\fin

\begin{remark}{\em  The general results of \cite{Ry} about stable measures entails that the exponential speed of convergence must be smaller than $\alpha/(1-\alpha)$ for arbitrary $\alpha$-stable processes with $0<\alpha<1$, and arbitrary semi-norms. The fact that in (\ref{ap1})
$$\frac{p\alpha}{p-\alpha}\; > \; \frac{\alpha}{1-\alpha}$$
as soon as $p <1$ is not contradictory with this result, since then $||.||_p$ is no more a semi-norm.}
\end{remark}

\subsection{The case when $|Z|$ is not a subordinator and $p > 1$.}

In this situation, the random variable $||Z||^p_p$ is more difficult to handle with. When $\alpha\, <\, 2\,\wedge\, p$ and in the symmetric case, Greenwood had proved that $||Z||^p_p$ belongs to the domain of attraction of some positive $(\alpha/p)$-stable law - see Corollary 2 in \cite{Gr}. In particular, the upper tails of $||Z||^p_p$ are described followingly: there exists $\k'_{\alpha, p}\in (0, +\infty)$ such that
$$\lim_{x\uparrow +\infty}x^{\alpha/p}\pb\lcr ||Z||_{p}\geq x \rcr
 \; =\; \k'_{\alpha,p}.$$
Greenwood had also raised the question whether $||Z||^p_p$ itself should be a positive $(\alpha/p)$-stable variable - see Section 6 in \cite{Gr}. If the answer to this question were yes, this would give an immediate answer to our small deviation problem, as in Proposition \ref{p21}. Unfortunately, the answer to this question is no:

\begin{theorem} \label{t22} Let $Z$ be a real stable process with index $\alpha\in (0,2]$. For every $p\, >\, 1\,\vee\,\alpha$, the positive random variable $||Z||^p_p$ is not stable.
\end{theorem}

Nevertheless, the main result of this paper states that the lower tails of $||Z||^p_p$ are similar to those of a positive stable random variable with index $\alpha/p$:

\begin{theorem} \label{t23} Let $Z$ be a real stable process with index $\alpha\in (0,2]$, such that $|Z|$ is not a subordinator. For every $p\, > \, 1\, \vee\,\alpha$, there exists $\k_{\alpha,p}\in (0, +\infty)$ such that
$$\lim_{\ee\downarrow 0}\ee^{\frac{p\alpha}{p-\alpha}}\log\pb\lcr ||Z||_{p}\leq \ee \rcr
 \; =\; -\k_{\alpha,p}.$$
\end{theorem}

Contrary to Proposition \ref{p21}, the small deviation constant is not explicit in Theorem \ref{t23}. However, thanks to (\ref{infvar}), we immediately see that
\begin{equation}\label{kap}
\k_{\alpha,p}\; \geq \; \lpa\frac{p-\alpha}{\alpha}\rpa {\lpa \lpa\frac{(c_- + c_+)}{p}\rpa\Gamma(1-\alpha/p)\rpa}^{\frac{p}{p-\alpha}}
\end{equation}
as soon as $\alpha < 2$. By continuity with Proposition 1, it seems natural to conjecture that the inequality is indeed an equality, but we were unable to prove this. On the other hand, the Gaussian case $\alpha = 2$ may well exhibit a different constant - see Paragraph 6.2. As a rule, the computation of small deviation constants has proved to be a hard task. Only few cases are known, almost always in a Gaussian framework (with the notable exception of symmetric stable processes under the $L_2$-norm, see Lemma 2.2 in \cite{Shi}). For Brownian Motion $W$ and $||.||_{1/p}$ the $(1/p)$-H\"older semi-norm on $[0,1]$ $(p >2)$, it was proved in \cite{BR} that
\begin{equation}\label{hoelder2}
\lim_{\ee\downarrow 0}\ee^{\frac{2p}{p-2}}\log\pb\lcr ||W||_{1/p}\leq \ee \rcr
 \; =\; - c_p\in (-\infty,0).
\end{equation}
The exact value of $c_p$ is still unknown, but some bounds are given in \cite{KL} with the help of Ciesielski's isomorphism theorem. It follows readily from (\ref{hoelder}) that 
$$\k(2, p)\; \leq\; c_p.$$
In this paper, an elementary use of Bochner's subordination allows us to get an upper bound on $\k_{\alpha,p}$ in the symmetric case, and when $\alpha < p< 2$. We note
$$D_{\alpha, p}\; =\; \lpa\frac{p-\alpha}{\alpha}\rpa {\lpa \lpa\frac{2\alpha \Gamma(\alpha) \sin (\pi\alpha/2)}{p\pi}\rpa\Gamma(1-\alpha/p)\rpa}^{\frac{p}{p-\alpha}}\kappa^{\frac{p}{p-\alpha}}$$
and
$$d_\alpha\; =\;\lpa\frac{2-\alpha}{2}\rpa \alpha^{\frac{\alpha}{2-\alpha}} \kappa^{\frac{2}{2-\alpha}}$$
for every $\alpha\in (0,2)$, where $\kappa$ is the parameter appearing in (\ref{kappa}).
\begin{theorem} \label{t24} Suppose that $Z$ is symmetric and that $\alpha < 2 < p$. Then
$$D_{\alpha, p}\; \leq \; \k_{\alpha,p}\; \leq \;c_p \; +\; d_{\alpha}.$$
\end{theorem}
The remainder of this paper is devoted to the proofs of Theorems 2, 3 and 4.

\section{Proof of Theorem \ref{t22}}

We first consider the case $\alpha < 2$. Recall that if $X$ is a positive $(\alpha/p)$-stable variable, then 
$$X\; \elaw\; 2^{-p/\alpha}(X_1 + X_2),$$
where $X_1$ and $X_2$ are independent copies of $X$. When $|Z|$ is a subordinator, it is clear from (\ref{pge1}) and the fact that $(a+b)^p > a^p + b^p$ for every $a, b >0$, that this latter equality fails. We now focus on the case when $|Z|$ is not a subordinator. For every $0\leq s\leq t$, we write
$${\VV}_p^{[s,t]}\; =\; \sup_{\tiny{t_0 < \ldots < t_k\in 
[s,t]}}\sum_{j=1}^k{\lva Z_{t_j} - Z_{t_{j-1}}\rva}^p.$$
Set $\VV_p = {\VV}_p^{[0,1]} = ||Z||^p_p$. By scaling, stationarity and independence of the increments of $Z$, it is clear that 
$2^{p/\alpha}{\VV}_p^{[0,1/2]}$ and $2^{p/\alpha}{\VV}_p^{[1/2,1]}$ are two independent copies of $\VV_p$. Hence, if $\VV_p$ were a positive $(\alpha/p)$-stable variable, we would have
$${\VV}_p^{[0,1/2]}\; +\; {\VV}_p^{[1/2, 1]}\;\elaw\; \VV_p.$$
Since obviously a.s. ${\VV}_p^{[0,1/2]}\; +\; {\VV}_p^{[1/2, 1]}\;\leq\; \VV_p$, this would entail
$${\VV}_p^{[0,1/2]}\; +\; {\VV}_p^{[1/2, 1]}\; =\; \VV_p\;\;\;\;\mbox{a.s.}$$
We will see that this latter equality is false. More precisely, we will prove that for every $0\leq s < t < u$ and every $K > 0$,
\begin{equation}
\label{false}
\pb\lcr {\VV}_p^{[s,u]}\, >\, K + {\VV}_p^{[s, t]} + {\VV}_p^{[t, u]}\rcr \; > \; 0.
\end{equation}
Indeed, since $p > 1$ and $s < t < u$, we can choose $M$ big enough and $\ee$ small enough such that
\begin{equation}
\label{meps}
{\lpa M(u-s) - \ee\rpa}^p - {\lpa M(t-s) + \ee\rpa}^p - {\lpa M(u-t) + \ee\rpa}^p\; >\; K.
\end{equation}
Consider the function $\phi : x \mapsto Z_s + M(x-s)$ from $[s,u]$ to $\R$. It is easy to see that
\begin{equation}
\label{norm}
{\lva\lva\phi\rva\rva}_{I,p}\; = \; M|I|
\end{equation}
for every interval $I\subset [s,u]$, where $|I|$ stands for the Lebesgue measure of $I$. On the other hand, the fact that $|Z|$ is not a subordinator entails that the (infinite) L\'evy measure $\nu$ is supported by the whole $\R$, and in particular by a neighbourhood of 0. So we can reason exactly as in Proposition 15 in \cite{pvaria} and get
\begin{equation}
\label{pb}
\pb\lcr \lva\lva Z - \phi\rva\rva_{[s,u], p} \, < \, \ee\rcr\; > \; 0.
\end{equation}
But it is clear that
\begin{eqnarray}
\label{triangle}
\lacc\lva\lva Z -\phi\rva\rva_{[s,u], p} \, < \, \ee\racc & \subset & \lacc\lva\lva Z -\phi\rva\rva_{[s,u], p} \, < \, \ee, \;\lva\lva Z - \phi\rva\rva_{[s,t], p} \, < \, \ee, \;\lva\lva Z - \phi\rva\rva_{[t,u], p} \, < \, \ee\racc\nonumber\\
           & \subset & \lacc \VV_p^{[s,u]} > {\lpa M(u-s) - \ee\rpa}^p\racc\;\cap\;\lacc\VV_p^{[s,t]} > {\lpa M(t-s) + \ee\rpa}^p\racc\nonumber\\
& & \;\;\;\;\;\;\;\;\;\;\;\;\;\;\;\;\;\;\;\;\cap\;\lacc\VV_p^{[t,u]} > {\lpa M(u-t) + \ee\rpa}^p\racc,
\end{eqnarray}
where in the second line we used (\ref{norm}) and triangle's inequality. Putting (\ref{meps}), (\ref{pb}) and (\ref{triangle}) together yields (\ref{false}), which entails that $\VV_p$ is not a positive $(\alpha/p)$-stable random variable. Since $\VV_p$ belongs to the domain of attraction of such a variable \cite{Gr}, we see that $\VV_p$ is not {\em any} stable variable at all. This completes the proof of Theorem \ref{t22} when $\alpha < 2$.

\vspace{2mm}

The case $\alpha = 2 < p$ is very easy, since then $||Z||_p$ is a finite Gaussian norm which must have some square exponential upper tails, so that $\VV_p = ||Z||^p_p$ cannot be a positive stable variable. This can also be viewed from the fact that if $\VV_p$ were stable, then $t\mapsto \VV^{[0,t]}_p$ would be a continuous increasing stable process, hence a deterministic linear function, which is obviously not the case. 
 
\fin

\begin{remark}{\em This result gives some insight on the random partition $\TT$ of $[0,1]$ realizing $\VV_p$. When $\alpha < 2$, it follows from Corollary 1 in \cite{Gr} that
$$\pb\lcr \mbox{$\TT$ is not dense in $[0,1]$}\rcr\; > \; 0.$$
It would be interesting to know if the latter probability is one. This is clearly true in the case $\alpha = 2 < p$, by continuity.}
\end{remark}

\section{Proof of Theorem \ref{t23}}

We use the same notations as in the preceding section concerning $\VV_p^{[s,t]}$ and $\VV_p$. We set $s_k^n\; = k/n$ and $\VV_p^{k,n}\; =\; \VV_p^{[s_{k-1}^n,s_k^n]}$ for every $n\geq 1$ and $k = 0,\ldots, n$. We will also consider
$$||Z||_{\omega}\; =\;\lim_{q\uparrow +\infty}||Z||_q\; =\; \sup_{0\leq s< t \leq 1}|Z_t- Z_s|,$$
see \cite{ChG} for the last inequality. In the remainder of this section, we will fix $p > 1$ once and for all.

\subsection{Three lemmas}

The proof of Theorem \ref{t23} relies on three easy lemmas. The first two are of deterministic nature, although we stated them with the process $Z$ for the sake of concision.

\begin{lemma}\label{omega} For every $n\geq 1$,
$$\VV_p\; \leq\; \sum_{k =1}^n \VV_p^{k,n}\; +\; n\, ||Z||_{\omega}^p\;\;\;\;\;\mbox{a.s.}$$
\end{lemma}

\noindent
{\em Proof.} We fix $n\geq 1$ and set $s_j = s^n_j$ for simplicity. Let $0 = t_0 < t_1 < \ldots < t_k = 1$ be a partition of $[0,1]$. Introduce $q^-_0 = q^+_0 = 0$ and set
$$q^-_j\;=\;\sup\lacc q\in\{0,\ldots, k\}\;/\; t_q <s_j \racc, \;\;\; q^+_j\;=\;\inf\lacc q\in\{0,\ldots, k\}\;/\; t_q\, \geq\, s_j\,\wedge\, 1\racc
$$ 
for every $j\geq 0$. Notice that $q^+_n = k = q^-_r = q^+_r$ for every $r>n$. Consider 
$$J\; = \; \lacc j\geq 0, \; | \;t_{q^+_{j+1}} > 
t_{q^+_j}\;\;\mbox{or}\;\; j =n\racc.$$
Notice that $\sharp\, J \,\leq\, n +1$ and that 
$$s_j \leq t_{q^+_j}\leq t_{q^-_{j+1}} < s_{j+1}$$
if $j\in J$. We get
\begin{eqnarray*}
\sum_{i = 1}^k {\lva Z_{t_i} - Z_{t_{i-1}}\rva}^p & = & \sum_{j\in J}\lpa\lpa\sum_{\tiny{t_{q^+_j}\leq t_q < t_{q+1} \leq 
t_{q^-_{j+1}}}}\!\!\!\!{\lva Z_{t_{q+1}} - Z_{t_{q}}  \rva}^p\rpa\; +\; {\lva Z_{t_{q^+_j}}- Z_{t_{q^-_j}}\rva}^p\rpa\\
& \leq & \sum_{k =1}^n \VV_p^{k,n}\; +\; n\, ||Z||_{\omega}^p\;\;\;\;\;\mbox{a.s.}
\end{eqnarray*}

\fin

\begin{lemma}\label{diff} For every $n\geq 1$ and $\ee > 0$,
$$\bigcap_{k=1}^n\lacc \VV_p^{k,n} \leq \ee,\; (Z_{s^n_k} - Z_{s^n_{k-1}})\, Z_{s^n_k} \leq 0\racc\; \subset\;\lacc \VV_p\leq 5^p n\ee\racc.$$

\end{lemma}

\noindent
{\em Proof.} Since $\lva Z_{s^n_k} - Z_{s^n_{k-1}}\rva^p \leq \VV_p^{k,n}$ for each $k = 1,\ldots, n$, we have
\begin{eqnarray*} \bigcap_{k=1}^n\lacc \VV_p^{k,n} \leq \ee,\; (Z_{s^n_k} - Z_{s^n_{k-1}})\, Z_{s^n_k} \leq 0\racc & \subset & \bigcap_{k=1}^n\lacc \VV_p^{k,n} \leq \ee,\; \lva Z_{s^n_k} \rva^p \leq \ee\racc\\
& \subset & \bigcap_{k=1}^n\lacc \VV_p^{k,n} \leq \ee\racc\; \cap\;\lacc ||Z||_{\omega}^p \leq 4^p \ee\racc \\
& \subset & \lacc \VV_p\leq 5^p n\ee\racc,
\end{eqnarray*}
where in the last inclusion we used Lemma \ref{omega}.

\fin

\begin{lemma}\label{positive} Suppose that $|Z|$ is not a subordinator. Then for every $x > 0$,
$$\pb\lcr ||Z||_p < x,\; Z_1 > 0\rcr\; > \; 0.$$

\end{lemma}

\noindent
{\em Proof.} Since $Z$ is strictly stable, it follows from Example 4 in \cite{pvaria} that for every $x > 0$,
$$\pb\lcr ||Z||_p < x \rcr\; > \; 0.$$
If $Z$ is symmetric, then clearly
$$\pb\lcr ||Z||_p < x,\; Z_1 > 0\rcr\; =\; \frac{1}{2}\lpa\pb\lcr ||Z||_p < x \rcr\rpa\; > \; 0,$$
since $||Z||_p = ||-Z||_p$ and $\pb[Z_1 = 0] = 0$. If $Z$ is not symmetric, we introduce the function $\phi_x : t\mapsto xt/2$ from $\R^+$ to $\R$. Since $|Z|$ is not a subordinator, its L\'evy measure is supported by the whole $\R$ and again, we can reason as in Proposition 15 of \cite{pvaria} to obtain
$$\pb\lcr ||Z-\phi_x||_p < x/4 \rcr\; > \; 0.$$
By triangle's inequality and the fact that $||\phi_x||_p = x/2$, we have
$$\lacc ||Z-\phi_x||_p < x/4\racc\; \subset\;\lacc ||Z||_p < x,\; Z_1 > x/4\racc.$$
This entails
$$\pb\lcr ||Z||_p < x,\; Z_1 > 0\rcr\; > \; 0$$
as desired.
\fin

\subsection{End of the proof}

We will first prove
\begin{equation}
\label{liminf}
\liminf_{\ee\downarrow 0}\ee^{\frac{p\alpha}{p-\alpha}}\;\pb\lcr ||Z||_p \le \ee \rcr\; > \; -\infty.
\end{equation}
When $\alpha =2$, this follows easily from Stolz's criterion - see Theorem 1 in \cite{St1}. When $\alpha <2$, we will use Lemmas \ref{diff} and \ref{positive}. Fix $\ee > 0$ and let $n\in\NN$ be such that 
$$\ee^{\frac{p\alpha}{\alpha -p}}\leq n \leq \ee^{\frac{p\alpha}{\alpha -p}} + 1.$$ 
By Lemma \ref{diff},
\begin{eqnarray*}
 \pb\lcr ||Z||_p \le \ee \rcr & \geq & \pb\lcr\bigcap_{k=1}^n\lacc\VV_p^{k,n} \leq c_0\,\ee^p/n,\; (Z_{s^n_k} - Z_{s^n_{k-1}})\, Z_{s^n_k} \leq 0\racc\rcr\\
& = & \pb\lcr\bigcap_{k=1}^n\lacc\VV_p^{k,n} \leq c_1\, n^{-p/\alpha},\; (Z_{s^n_k} - Z_{s^n_{k-1}})\, Z_{s^n_k} \leq 0\racc\rcr
\end{eqnarray*}
where $c_0, c_1$ are positive finite constants not depending on $n$. Define
$$\AA_j \; =\; \bigcap_{k=1}^j\lacc \VV_p^{k,n} \leq c_1 n^{-p/\alpha},\; (Z_{s^n_k} - Z_{s^n_{k-1}})\, Z_{s^n_k} \leq 0\racc$$ 
and $\FF_j$ the $\sigma$-algebra generated by $\lacc Z_t, \; t\leq s^n_j\racc$, for every $j = 1,\ldots, n$. On $\lacc Z_{s^n_{n-1}} \leq 0\racc$, we see that
\begin{eqnarray*}
 \esp\lcr \AA_n\, |\, \FF_{n-1}\rcr & = & {\Un}^{}_{\AA_{n-1}}\pb\lcr\VV_p^{n,n} \leq c_1 n^{-p/\alpha},\; Z_{s^n_n} - Z_{s^n_{n-1}}\geq 0 \rcr\\
& = & {\Un}^{}_{\AA_{n-1}}\pb\lcr\VV_p\leq c_1,\; Z_1\geq 0 \rcr
\end{eqnarray*}
by scaling, independence and stationarity of the increments of $Z$. Similarly, we have
\begin{eqnarray*}
 \esp\lcr {\Un}^{}_{\AA_n}\, |\, \FF_{n-1}\rcr & = & {\Un}^{}_{\AA_{n-1}}\pb\lcr\VV_p\leq c_1,\; Z_1\leq 0 \rcr
\end{eqnarray*}
on $\lacc Z_{s^n_{n-1}} \geq 0\racc$. By a direct induction,
$$\pb\lcr ||Z||_p \le \ee \rcr\; = \; \pb\lcr \AA_n\rcr\; \geq \;c_2\, \pb \lcr\AA_{n-1}\rcr\; \geq \;c_2^n,$$ 
where we set
$$c_2\; = \; \min\lacc\pb\lcr\VV_p\leq c_1,\; Z_1\leq 0 \rcr, \pb\lcr\VV_p\leq c_1,\; Z_1\geq 0 \rcr\racc,$$
which is positive by Lemma \ref{positive}. Hence
$$\ee^{\frac{p\alpha}{p-\alpha}}\log \pb\lcr ||Z||_p \le \ee \rcr\; \geq \; n \ee^{\frac{p\alpha}{p-\alpha}}\log c_2 \; \geq \; \log c_2 \; > \; -\infty$$
for every $\ee > 0$, and (\ref{liminf}) is proved.

\vspace{2mm}

We now proceed to the last step of the proof. If $Z$ is symmetric, it follows from Theorem 3.1 in \cite{LiSi} that
$$\lim_{\ee\downarrow 0}\ee^{\frac{p\alpha}{p-\alpha}}\;\pb\lcr ||Z||_p \le \ee \rcr\;$$
exists and belongs to $[-\infty,0)$. But we know from (\ref{liminf}) that it cannot be $-\infty$, and the proof is complete. 

If $Z$ is not symmetric, we can actually repeat almost verbatim the proof of Theorem 3.1 in \cite{LiSi}, to obtain the existence of the constant. We give the details for completeness. Set $q = p/\alpha > 1$. Since $p > 1$, it is clear that 
$$\VV_p\; \geq\; {\VV}_p^{[0,t]}\; +\; {\VV}_p^{[t, 1]}\;\;\;\mbox{a.s.}$$
for every $t\in (0,1)$. On the other hand, $t^{-p/\alpha}{\VV}_p^{[0,t]}$ and $(1-t)^{-p/\alpha}{\VV}_p^{[1-t,1]}$ are two independent copies of $\VV_p$, and we easily see that
$$\esp\lcr e^{- (x+y)^q \VV_p} \rcr \;\le\; \esp\lcr e^{-x^q
\VV_p} \rcr \ \esp\lcr e^{-y^q \VV_p} \rcr $$
for every $x, y \geq 0$. If we now set
$$\Psi(h) = \log \esp \lcr e^{-h^q \VV_p} \rcr$$
for every $h \geq 0$, this entails that $\Psi$ is a continuous negative function which satisfies
$\Psi(x+y) \leq \Psi(x) + \Psi(y)$ for every $x, y \geq 0$. By the standard subadditivity argument, we obtain
$$\lim_{h\to\infty} \frac{\Psi(h)}h \; = \; \inf_{h \geq 0}
\frac{\Psi(h)}h\; = \; -C \;\in \; [-\infty, 0),$$
which reads
$$\lim_{\lambda\to +\infty}\lambda^{1/q}\log \esp\lcr e^{-\lambda
\VV_p} \rcr\; =\; -C.$$
Since $q>1$, de Bruijn's exponential Tauberian theorem entails that
$$\lim_{\ee\downarrow 0}\ee^{\frac{p\alpha}{p-\alpha}}\log\pb\lcr ||Z||_p \le \ee \rcr$$
exists, and is finite because of (\ref{liminf}). This completes the proof.

\fin

\section{Proof of Theorem \ref{t24}}

We fix $p > 2$. The first inequality $\k_{\alpha,p} \geq D_{\alpha, p}$ follows readily from (\ref{kap}) and (\ref{cplus}). To prove that
\begin{equation} \label{cpda}
\k_{\alpha,p} \leq c_p + d_{\alpha},
\end{equation}
we first remark that since $Z$ is symmetric, it can be rewritten $Z = W\,\circ\,\sigma$, where $\lacc W_t,\; t\geq 0\racc$ is a standard linear Brownian motion and $\lacc \sigma_t,\; t\geq 0\racc$ an independent $(\alpha/2)$-stable subordinator, whose Laplace transform is given by
\begin{equation}\label{kappa2}
\esp\lcr e^{-\lambda Z_t}\rcr\; =\; e^{-\kappa 2^{\alpha/2} t\lambda^{\alpha/2}}.
\end{equation}
For every $T > 0$ we introduce
$$\HH(p,T)\; =\;\sup_{0\leq s< t \leq T}\frac{|W_t-W_s|}{{|t-s|}^{1/p}},$$
which is a.s. finite since $p > 2$. Let $0=t_0 < \ldots < t_n =1$ be any partition of $[0,1]$. We have
\begin{eqnarray*}
\sum_{i=1}^n \lva Z_{t_i} - Z_{t_{i-1}}\rva^p & = & \sum_{i=1}^n \lva W_{\sigma_{t_i}} - W_{\sigma_{t_{i-1}}}\rva^p\\
& \leq & \HH(p,\sigma_1)^p \sum_{i=1}^n \lva \sigma_{t_i} - \sigma_{t_{i-1}}\rva\\
& \elaw & \HH(p,1)^p \sigma_1^{p/2 - 1}\sum_{i=1}^n \lva \sigma_{t_i} - \sigma_{t_{i-1}}\rva\\
& = & \HH(p,1)^p \sigma_1^{p/2},
\end{eqnarray*}
where in the third line we used the scaling property of $W$ and the fact that $W$ and $\sigma$ are independent. In particular
$$\pb \lcr \VV_p \leq \ee \rcr\;\geq\;\pb\lcr  \HH(p,1)^p \sigma_1^{p/2}\leq\ee\rcr$$
for every $\ee > 0$. Set
$$r \; =\; \frac{\alpha(p-2)}{2(p-\alpha)}\; \in\; (0,1),\;\;\;\; 1 - r \; =\; \frac{p(2-\alpha)}{2(p-\alpha)},$$
and $\HH_p = \HH(p,1)$ for concision. By independence, it is clear that for every $\ee > 0$,
$$\pb\lcr\VV_p \le \ee \rcr\; \geq \;\pb\lcr\HH^p_p \le \ee^r \rcr\pb\lcr\sigma_1^{p/2} \le \ee^{1-r}\rcr.$$
Hence,
\begin{equation}\label{ineg}
\ee^{\frac{p\alpha}{p-\alpha}}\log\pb\lcr ||Z||_p \le \ee \rcr \; \geq \; \ee^{\frac{\alpha}{p-\alpha}}\log\pb\lcr \HH_p \le \ee^{r/p}\rcr\; + \;\ee^{\frac{\alpha}{p-\alpha}}\log\pb\lcr \sigma_1 \le \ee^{2(1-r)/p}\rcr.
\end{equation}
On the one hand, setting $x = \ee^{r/p}$, it follows from (\ref{hoelder2}) that
\begin{eqnarray}\label{hp}
\ee^{\frac{\alpha}{p-\alpha}}\log\pb\lcr \HH_p \le \ee^{r/p}\rcr & = & x^{\frac{2p}{p-2}}\log\pb\lcr \HH_p \le x \rcr\;\rightarrow \; - c_p\;\;\;\mbox{as $\ee\downarrow 0$.}
\end{eqnarray}
On the other hand, setting $x = \ee^{2(1-r)/p}$, (\ref{cplus}) (\ref{subo}) and (\ref{kappa2}) entail that
\begin{eqnarray}\label{sigma}
\ee^{\frac{\alpha}{p-\alpha}}\log\pb\lcr \sigma_1 \le \ee^{(2-\alpha)/(p-\alpha)}\rcr & = & x^{\frac{\alpha}{2-\alpha}}\log\pb\lcr \sigma_1 \le x \rcr \; \rightarrow - d_{\alpha}\;\;\;\mbox{as $\ee\downarrow 0$.}
\end{eqnarray}
Putting (\ref{ineg}), (\ref{hp}) and (\ref{sigma}) together yields (\ref{cpda}) as desired.

\fin

\section{Final remarks}

\subsection{General stable processes}

General stable processes can be viewed as drifted strictly stable processes (see e.g. Chapter 3 in \cite{sato} for an explanation of the importance of these processes in terms of limit distributions). If $\lacc W^{\mu}_t, \; t\geq 0\racc$ is a standard Brownian motion with drift $\mu\in\R$, it follows readily from the Cameron-Martin formula that
$$\lim_{\ee\downarrow 0}\ee^{\frac{p\alpha}{p-\alpha}}\log\pb\lcr ||W^{\mu}||_p \le \ee \rcr\; =\;-\k_{2,p}.$$
We were unable to obtain the same result in the Non-Gaussian framework. First, the path-transformation $Z \mapsto Z^{\mu} = \lacc Z_t + \mu \, t, \; t\geq 0\racc$ is no more absolutely continuous - see Theorem 33.1 in \cite{sato}. Second, the lack of self-similarity for $Z^{\mu}$ is quite cumbersome. For example, if ${\VV}_p^1$ and ${\VV}_p^2$ are independent and distributed like $\VV_p$, it is no more true that 
$$\VV_p\; \geq\; t^{-p/\alpha}{\VV}_p^1\; +\; {(1-t)}^{-p/\alpha}{\VV}_p^2$$
in law for every $t\in (0,1)$, which is a crucial step in proving the existence of the small ball constant.

\subsection{More on the small ball constant in the symmetric case}

In this paragraph we fix $p\, > \, 1\,\vee\, \alpha$ and we suppose that $Z$ is symmetric. Every identification or inequality will be made almost surely, and we skip the notation a.s. for concision. Consider the decreasing random sequence $\lacc\WW_{p,n}, \; n\geq 0\racc$ defined by
$$\WW_{p,n}\; = \;\sup_{\tiny{\begin{array}{c}
        t_0 < \ldots < t_k\in I\\
        |t_j - t_{j-1}|\leq 1/n
        \end{array}}}\lpa\sum_{j=1}^k{\lva Z_{t_j} - Z_{t_{j-1}}\rva}^p\rpa.$$
When $\alpha <2$, Greenwood had showed (see Corollaries 1 and 2 in \cite{Gr}) that $\VV_p$ belongs to the domain of attraction of the $(\alpha/p)$-stable variable
\begin{equation}\label{wpn}
\WW_p\; = \; \lim_{n\uparrow +\infty}\WW_{p,n}.
\end{equation}
We first remark that
\begin{equation}\label{ws}
\WW_p \; =\; S^p_1
\end{equation}
with the notations of Section 2. Indeed, set $Z^{\eta}$ for the compound Poisson process obtained from $Z$ in removing its jumps with size smaller than $\eta$, for every $\eta > 0$. If we define $\WW^{\eta}_p$ and $S^{p,\eta}_1$ similarly as above with respect to $Z^{\eta}$, then it is easy to see that $\WW^{\eta}_p \; =\; S^{p,\eta}_1$ for every $\eta >0$, since by symmetry $Z^{\eta}$ is a step-function with a finite number of steps. It is well-known and easy to see that $\lim_{\eta\downarrow 0} S^{p,\eta}_1\; = \; S^p_1$. Besides, the fact that $\WW^{\eta}_p \; \rightarrow \;\WW_p$ along a subsequence $\{\eta\}$ follows immediately from Theorem II in \cite{Bre}, so that (\ref{ws}) holds. 

\vspace{2mm}

Notice that for every $n \geq 1$,
$$\VV_p\; \geq\;\WW_{p,n}\; \geq\; \VV_p^{1,n}\, + \ldots + \, \VV_p^{n,n}\; \elaw\; n^{-p/\alpha}\lpa \VV_p^1\, + \ldots + \, \VV_p^n\rpa,$$
where in the last equality $\VV_p^1,\ldots,\VV_p^n$ stand for $n$ independent copies of $\VV_p$. De Bruijn's exponential Tauberian theorem entails
$$\lim_{\ee\downarrow 0}\ee^{\frac{\alpha}{p-\alpha}}\log\pb\lcr \VV_p \le \ee \rcr\; =\; \lim_{\ee\downarrow 0}\ee^{\frac{\alpha}{p-\alpha}}\log\pb\lcr n^{-p/\alpha}\lpa \VV_p^1\, + \ldots + \, \VV_p^n\rpa\le \ee \rcr,$$
which yields
\begin{equation}
\label{vwpn}
\lim_{\ee\downarrow 0}\ee^{\frac{\alpha}{p-\alpha}}\log\pb\lcr \VV_p \le \ee \rcr\; =\; \lim_{\ee\downarrow 0}\ee^{\frac{\alpha}{p-\alpha}}\log\pb\lcr \WW_{p,n} \le \ee \rcr
\end{equation}
for every $n\geq 1$. On the one hand, it follows from (\ref{wpn}) and (\ref{ws}) that
$$D_{\alpha, p}\; =\; \lim_{\ee\downarrow 0}\lpa\lim_{n\uparrow \infty}\ee^{\frac{\alpha}{p-\alpha}}\log\pb\lcr \WW_{p,n} \le \ee \rcr\rpa.$$
On the other hand, (\ref{vwpn}) obviously entails
$$\k_{\alpha, p}\; =\; \lim_{n\uparrow \infty}\lpa\lim_{\ee\downarrow 0}\ee^{\frac{\alpha}{p-\alpha}}\log\pb\lcr \WW_{p,n} \le \ee \rcr\rpa.$$ 
Twisting the limits in $\ee$ and $n$ would lead to the conjectured equality
$$\k_{\alpha, p}\; =\; D_{\alpha,p}$$ 
when $\alpha < 2$, but we were unable to decide whether this is plausible or not. Notice that this is clearly false when $\alpha = 2$, since then $\WW_{p,n} \rightarrow 0$ as $n\uparrow \infty$, whereas $\k_{2, p} > 0$.

\subsection{Some other two-sided estimates in the symmetric case}

In the symmetric case, the method used to prove Theorem \ref{t24} also applies to the supremum and oscillation semi-norms, and yields interesting upper estimates on the respective small ball constants which seem unnoticed in the literature. We use the notations
$$\LL_2(Z)\; =\; {\lpa\int_0^1 Z_t^2\, dt\rpa}^{1/2},\;\;\; ||Z||_{\infty}\; =\; \sup_{0\leq t\leq 1}|Z_t|\;\;\;\mbox{and}\;\;\; ||Z||_{\omega}\; =\; \sup_{0\leq s < t\leq 1}|Z_t -Z_s|.$$
Notice that since $Z$ starts from 0,
\begin{equation}\label{l2}
\LL_2 (Z)\;\leq\;||Z||_{\infty}\;\leq\;||Z||_{\omega}\;\;\;\;\mbox{a.s.}
\end{equation}
Recently, it was proved by Shi (see Lemma 2.2. in \cite{Shi}) that 
\begin{equation}\label{l22}
\lim_{\ee\downarrow 0}\ee^{\alpha}\log\pb\lcr \LL_2 (Z)\le \ee\rcr\; =\; -\gamma_{\alpha}\; \dde\; -2^{\alpha/2}\kappa \lpa\inf_{f\in\AAA}\int_{\R} |x|^{\alpha}f(x)\, dx\rpa,
\end{equation}
where $\AAA$ stands for the set of absolutely continuous probability densities $f$ on $\R$ satisfying
$$\frac{1}{8}\lpa\int_{\R} \frac{(f'(x))^2}{f(x)}\, dx\rpa\; \leq \; 1.$$
It is well-known that $\gamma_2 = \kappa/4$. A long time ago, it was proved by Taylor \cite{Ta} and Mogul'ski\v {\i} \cite{Mo} that
$$\lim_{\ee\downarrow 0}\ee^{\alpha}\log\pb\lcr ||Z||_{\infty} \le \ee\rcr\; =\; -\k_{\alpha, \infty}\; \in\; (-\infty, 0).$$
It is very well-known that $\k_{2, \infty} = \pi^2/8$, but no formula for $\k_{\alpha, \infty}$ seems to be available when $\alpha < 2$. Last, it follows readily from the standard subadditivity argument that
$$\lim_{\ee\downarrow 0}\ee^{\alpha}\log\pb\lcr ||Z||_{\omega} \le \ee\rcr\; =\; -\k_{\alpha, \omega}\; \in\; [-\infty, 0).$$
The distribution of $||Z||_{\omega}$ was computed by Feller \cite{Fe} in the Gaussian case and yields $\k_{2, \omega} = \pi^2/2$, but again, no formula for $\k_{\alpha, \omega}$ is available when $\alpha < 2$.  

\begin{prop} The following inequalities hold:
$$\gamma_{\alpha}\; \leq\; \k_{\alpha, \infty}\;\leq\;\pi^2/8\, +\,d_{\alpha}\;\;\;\mbox{and}\;\;\;\k_{\alpha, \infty}\; \leq\; \k_{\alpha, \omega}\;\leq\;\pi^2/2\, +\,d_{\alpha}.$$
\end{prop}

\noindent
{\em Proof.} It follows readily from (\ref{l2}) and (\ref{l22}) that
$$\gamma_{\alpha}\; \leq\;\k_{\alpha, \infty}\; \leq\; \k_{\alpha, \omega}.$$
Using the same notations as in the proof of Theorem \ref{t24} and Brownian scaling, it is easy to see that for every $\ee > 0$ and every $r\in (0,1)$
$$\pb\lcr  ||Z||_{\infty} \le \ee \rcr\; \geq \;\pb\lcr ||W||_{\infty}\le \ee^r \rcr\pb\lcr\sigma_1\le \ee^{2(1-r)}\rcr.$$
In particular, setting $r =\alpha/2$, $x = \ee^{\alpha/2}$ and $y = \ee^{2-\alpha}$, we get
$$\ee^{\alpha}\log\pb\lcr ||Z||_{\infty}\le \ee \rcr\; \geq \; x^2\log\pb\lcr ||W||_{\infty}\le x\rcr\; + \;y^{\frac{\alpha}{2-\alpha}}\log\pb\lcr \sigma_1 \le y\rcr.$$
Letting $\ee$ tend to 0 yields
$$\k_{\alpha, \infty}\; \leq\; \pi^2/8 \; +\; d_{\alpha}.$$
The last inequality 
$$\k_{\alpha, \omega}\; \leq\; \pi^2/2 \; +\; d_{\alpha}$$
follows in the same way, replacing $||Z||_{\infty}$ and $||W||_{\infty}$ by $||Z||_{\omega}$ and $||W||_{\omega}$, respectively. 

\fin

\begin{remark}{\em The following small ball estimate, which is probably long well-known, was obtained in \cite{LiSi} (see Section 6 therein):
$$\lim_{\ee\downarrow 0}\ee^{\alpha}\log\pb\lcr \LL_p(Z) \le \ee\rcr\; =\; -{\tilde \k}_{\alpha, p}\; \in\; (-\infty, 0)$$
for every $p\geq 1$, with the obvious notation for $\LL_p(Z)$. However, the above subordination method does not yield any further information on ${\tilde \k}_{\alpha, p}$. We just have the obvious inequality
$${\tilde \k}_{\alpha, p}\; \leq \; \k_{\alpha, \infty}.$$
}
\end{remark}

\subsection{A lower estimate on the H\"older small ball constants for Brownian motion}

As an immediate consequence of Theorem 4, we get
$$c_p \; \geq \; \max_{\tiny{\begin{array}{c} \alpha \in (0, 2)\\ \kappa > 0\end{array}}}\lacc D_{\alpha, p}(\kappa) - d_{\alpha}(\kappa)\racc$$
(here we also take into account the fact that $D_{\alpha, p}$ and $d_{\alpha}$ depend on a scaling parameter $\kappa > 0$). Using the formula
$$\max_{x > 0} \lacc ax^q - bx^r\racc\; =\; \lpa 1/q - 1/r\rpa {(qa)}^{\frac{r}{r-q}}{(rb)}^{\frac{q}{q-r}}$$
which is valid for every $1 < q < r$ and $a, b > 0$, elementary (but tedious) calculations lead to
$$c_p \; \geq \;c_p' \; \dde \;\lpa\frac{p-2}{p}\rpa p^{\frac{2}{2-p}}\, 2^{\frac{p+2}{p-2}}\lpa\max_{\tiny{\alpha\in (0,2)}}{\lpa \frac{\Gamma\lpa\frac{1+\alpha}{2}\rpa\Gamma\lpa 1 - \frac{\alpha}{p}\rpa}{\Gamma \lpa \frac{1}{2}\rpa\Gamma\lpa 1 - \frac{\alpha}{2}\rpa}\rpa}^{\frac{2p}{\alpha (p-2)}}\rpa.$$
It is somewhat surprising that this valuable information on the H\"older small ball constants for Brownian motion can be obtained from considerations on Non-Gaussian stable processes. Unfortunately, one can show that
$$c_p'\; < \; \frac{1}{4}\, \Gamma\lpa\frac{p}{p-2}\rpa\; < \; 2^{\frac{2(1-p)}{p-2}}\int_0^{\infty}\lpa\frac{u^{2p/p-2}e^{-u^2/2}}{\int_0^u e^{-v^2/2}\, dv}\rpa du,$$
so that our lower bound is smaller than the one obtained in \cite{KL}. 
 
\vspace{2mm}

\noindent
{\em Acknowledgement:} Part of this work was done during a nice stay at the Sobolev Institute of Mathematics in Akademgorodok (Russia). I thank Anatoli\v {\i} Mogul'ski\v {\i} for his warm hospitality and excellent working conditions.

\end{document}